\documentclass[11pt,leqeq]{article}
\usepackage{amssymb,amsthm}
\usepackage{amsmath}
\usepackage{amscd}
\usepackage{xy}
\xyoption{all}
\usepackage[dvips]{graphicx}
\newtheorem{thm}{Theorem}

\newtheorem{prop}[thm]{Proposition}

\newtheorem{exmp}[thm]{Example}
\date{}

\begin{document}
\setlength{\baselineskip}{16pt}
\title{Indirect   Influences}
\author{Rafael D\'iaz} \maketitle

\begin{abstract}
We introduce the PWP method for counting indirect influences, and
compare it with three well-known methods, namely, MICMAC, Heat Kernel, and PageRank.
We provide combinatorial as well as probabilistic interpretation for the PWP
method.

\end{abstract}

\section{Introduction}

Our goal in this note is to compare four alternative approaches
to count indirect influences.  Two of these approaches,  the
MICMAC and PageRank methods, are well-tested and each
possesses a host of real-life applications. The Heat Kernel method has found interesting mathematical applications.
We believe that the PWP method,
while still in its infancy, may also be useful in a
variety of contexts. Our guide for the definition of the PWP method is threefold:
\begin{itemize}
\item Direct influences generate, by concatenation, indirect influences.

\item All indirect influences arise from the concatenation of direct influences.

\item Longer concatenations exert a lesser indirect influence.
\end{itemize}

We consider a discrete scenario where the variables involved are indexed, and often identified, with the set $[n] = \{1,2,...,n \}$ for some $n \in \mathbb{N}_+.$
We assume that direct influences among our variables are encoded into a
$\mathbb{R}$-weighted directed graph, without multiple edges, with $[n]$ as its
set of vertices. Thus a direct influence from the variable $j \in [n]$ on
the variable $i\in [n]$ corresponds with an edge of weight $D_{ij}$ from vertex $j$ to
vertex $i$. As usual the graph is encoded into the matrix  of
direct influences $D$, which is such that $D_{ij}$ is the weight of
the edge going from $j$ to $i$, if there is one, and $D_{ij}=
0$ otherwise. Our problem consists in evaluating the indirect
influences between the variables, i.e. finding the  matrix $T$ such that $T_{ij}$ measures the
indirect influence of variable $j$ on variable $i$. In this work we make a comparative analysis
of four possible definitions for $T$. \\

\section{MICMAC}

The MICMAC method, introduced by Godet \cite{godet},
works as follows. Let $D$ be the matrix associated
with the graph of direct influences. The graph of indirect influences is represented
by the matrix $T = D^k,$ where $k$ is a fixed small natural number, say $4$
or $5$.  The vector of indirect dependencies $d=(d_1,....d_n)$ and
the vector indirect of influences  $f=(f_1,....f_n)$, are such that their coefficients $d_j$ and $f_j$ are given, respectively, by:
\begin{equation}\label{df}
d_j = \sum_{i=1}^n T_{ji} \ \ \ \ \ \ \ \mbox{ and } \ \ \ \ \ \ \ f_j = \sum_{i=1}^n T_{ij}.
\end{equation}
Thus $(d_j, f_j)$ is equal to the $(\mathrm{in}, \mathrm{out})$-degree of the vertex
$j$ in the graph of indirect influences. The numbers $d_j$ and
$f_j$ encode valuable information,  for example, the most
influential variable is the one for which $f_j$ reaches its highest value.

\section{PageRank}

The celebrated method PageRank \cite{page} for computing indirect dependencies has already
revolutionized the world of internet search engines trough the remarkable success of Google.
The mathematics encoding the basic structure behind PageRank are however surprisingly simple. With hindsight we observe that PageRank differs from MICMAC in three main
ideas:
\begin{itemize}

\item Normalization of influences, i.e. use of column stochastic matrices.

\item Use of complete graphs, i.e. graphs whose associated matrice have no vanishing entries.

\item Taking infinite potencies of matrices.

\end{itemize}

Let  $D$ be the $n \times n$ matrix of direct influences whose entries are non-negative real numbers. Assume that the sum of the entries of each column of $D$ is either $0$ or $1$. With PageRank the matrix of indirect influences $T=T(D)$  is computed as follows:
\begin{equation}\label{pr1}
T = T(D)= \lim_{k \to \infty}\left[p\overline{D} + (1-p)E_n\right]^k
\end{equation}
where:
\begin{itemize}
\item  The parameter $ 0 < p <1$ is a chosen number close to $1$, say $p=0.86.$

\item The matrix $\overline{D}$ is obtained from $D$ by replacing the entries of each zero column of $D$
by $\frac{1}{n}$.

\item The matrix $E_n$ has all of its entries equal to $\frac{1}{n}$.
\end{itemize}

For the web application the matrix $D$ is constructed as follows. Consider
the web graph whose set of vertices $[n]$ represents the set of
all web pages in the world wide web. There is an edge from page
$j$ to page $i$ in the web graph if an hyperlink to page $i$
appears in the content of page $j.$ The matrix $D$ of direct
influences is given by:

$$D_{ij} =
  \left\lbrace
  \begin{array}{l}
     \frac{1}{\mathrm{out}(j)} \ \ \text{ if there is an edge from} \ j \ \text{to} \ i, \\
   \ \  0 \ \ \ \ \ \ \  \text{if there is no edge from} j \ \text{to} \ i.
 \\
  \end{array}
  \right.$$

Since $p \overline{D} + (1-p)E_n $ is a column stochastic matrix, then so is $T$ and thus the
vector of influences for $T$ is $(1,...,1)$. The entries of
the vector  $\frac{d}{n}=\frac{1}{n}(d_1, ... , d_n )$, where $d$ is the vector of dependencies of $T$ whose coordinates
are given by (\ref{df}), yields the PageRank
number of each web page. The greater the PageRank number of a page
the greater its importance. By construction  $p \overline{D} +
(1-p)E_n$ is the transition matrix of a primitive irreducible
Markov chain  \cite{meyer}.  Therefore the vector of indirect
dependencies $d$ is an eigenvector of $T$ with
eigenvalue $1$, i.e.  we have that $$Td^t=d^t.$$ Moreover, the matrix $T$ can be computed in terms of $d$ as follows
$$T= \left( \frac{d^t}{n}  \ \ \frac{d^t}{n} \ \ .... \ \ \frac{d^t}{n} \right).$$

\section{Heat Kernel}

The Heat Kernel method of Chung \cite{chung} proceeds as follows. Just as MICMAC depends on a parameter $k$, and PageRank on the parameter $p$,
the Heat Kernel method depends on a parameter $\lambda $, a fixed positive real number.
Given $D$, the matrix of direct influences, the
matrix $T=T(D)$ of indirect influences is given for a fixed parameter $\lambda >0$ by
$$T=T(D)= e^{\lambda(D - I)} ,$$ where $I$ is the identity matrix of the appropriate size.
Note that when $D=0$, i.e. in the absence of direct influences, the matrix of indirect influences $T=T(0)$ is not the zero matrix, since each vertex self-influences itself. Indeed, in this case we have that
$$T=T(0)= e^{-\lambda}I.$$

\section{PWP}

Let us introduce the PWP  method for counting indirect influences. Unlike PageRank, the PWP method
can be applied to any matrix of direct influences, even matrices with negatives entries. Theorem \ref{ja}
below shows that whereas MICMAC focuses on paths of a fixed length $k$, and PageRank
focuses on infinite long paths, the Heat Kernel and PWP methods take into account paths of various lengths. The PWP method
avoids the self-influences that are included, by default, in the Heat Kernel method.\\

In a nutshell the PWP method can described as follows. As with the Heat Kernel, we fixed a parameter
$\lambda >0$. Assume we are given the matrix $D$ of direct
influences, and let $T$, the matrix of indirect influences, be given by
$$T = T(D) = \frac{e^{\lambda D} - I}{e^{\lambda} - 1} = \frac{e_+^{\lambda D}}{e_+^{\lambda}},$$
 where $$e_+^x= e^x -1=\sum_{k=1}^{\infty}\frac{x^k}{k!}.$$

Note that $D^k= \frac{\partial^k}{\partial \lambda^k}\left(
e_+^{\lambda}T \right)|_{\lambda =0}$, for $k \geq 1$, and thus in principle one
can compute the MICMAC matrix of indirect influences from the PWP
matrix of indirect influences. \\

Let us consider a rather trivial example which however highlights some of the
differences between the four methods.
Let $\mathrm{L}_2$ be the graph $1 \rightarrow 2, \ $ i.e. the graph
with vertex set $[2]$ and a unique edge from $1$ to $2$. \\

The MICMAC matrix of indirect influences in this case vanishes for
$k \geq 2$, thus no vertex influences or depends on another vertex.\\

PageRank makes $1$ the most influential vertex, and $2$ the most dependent vertex. Note the difference
with MICMAC. It assigns a non-vanishing dependency of vector $1$ on vector $2$. It also assigns
a non-vanishing dependency of $1$, respectively $2$, on itself.\\

Heat Kernel makes $1$ the most influential vertex, and $2$ the most dependent vertex.
It assigns a vanishing dependency of $1$ on $2$. Note the difference
with PageRank. It assigns a non-vanishing dependency of $1$, respectively $2$, on itself.\\

With PWP the only non-vanishing entry of the matrix $T$ is  $$T_{12}=\frac{\lambda}{e^{\lambda} - 1},$$ thus making $1$ the most
influential vertex, and $2$ the most dependent vertex. Note the differencie with MICMAC. Vertex $2$ does not exert any
influence on vertex $1$, unlike with PageRank. There is not self-influence of a vertex on itself,
which illustrates  the difference between the PWP and Heat Kernel methods.\\

Note that the indirect influence of vertex $1$ over vertex $2$ is given by
$$\frac{\lambda}{e^{\lambda} - 1} = \sum_{k=0}^{\infty}B_k\frac{\lambda^k}{k!},$$
where the coefficients $B_k$ are the so-called Bernoulli numbers.
See \cite{blan} for an explicit definition and a combinatorial
interpretation for the Bernoulli numbers. Note also that
$$\frac{\lambda}{e^{\lambda} - 1}\leq 1, \ \ \ \ \ \ \ \frac{\lambda}{e^{\lambda} - 1}\rightarrow 1 \ \ \ \mbox{as} \ \ \ \lambda \rightarrow 0,
\ \ \ \ \ \ \text{and} \ \ \ \ \ \frac{\lambda}{e^{\lambda} - 1}\rightarrow 0 \ \ \ \text{as} \ \ \ \lambda \rightarrow \infty.$$  Therefore, the indirect influence that
$1$ exerts over $2$ is lesser than the original direct influence,
it approaches its original value when $\lambda$ approaches $0$,
and it is negligible when $\lambda$ is a large number.  In other words, the
direct influence that $1$ exerts over $2$  becomes less relevant
as $\lambda$ increases, since with the PWP method the
paths of length close to $\lambda$ are the most relevants.\\

Next result states some of the basic properties of the map $D
\rightarrow T(D)$ and provides a probabilistic interpretation
for the PWP method.

\begin{thm}{\em Let $T=T(D)$ be the matrix of indirect influences with the PWP method.
\begin{enumerate}

\item If $D=0$, then $T(D)=0.$ \ \ \  $T(D^t) = T(D)^t.$

\item $T(QCQ^{-1})=QT(C)Q^{-1}$.

\item Let $D_1 \oplus D_2 \in M_m(\mathbb{R})\oplus
M_n(\mathbb{R})\subseteq  M_{m+n}(\mathbb{R}) $, then
$T(D_1 \oplus D_2)=T(D_1) \oplus T(D_2)$.

\item If $D$ is a column stochastic matrix, then $T(D)$ is a column stochastic matrix.

\item If $[D_1, D_2] = 0$, then we have that $$T(D_1 + D_2)= e_+^{\lambda}T(D_1)T(D_2) + T(D_1) + T(D_2).$$

\item Let $D \in M_n(\mathbb{R})$, then $T(D)$ is the expected matrix
of the random matrix $\widehat{D},$ where
\begin{enumerate}
\item $\widehat{D}: \mathbb{N}_+ \longrightarrow
M_n(\mathbb{R})$  is the random matrix
given by $\widehat{D}(k)=D^k$.

\item $\mathbb{N}_+$ is the probability space with probability function   $p(k)=\frac{\lambda^k}{e_+^{\lambda}k!}$.

\end{enumerate}

\end{enumerate}
}
\end{thm}

\begin{proof}
1. $ T(D^t) = \frac{e_+^{\lambda D^t}}{e_+^{\lambda}} = \left( \frac{e_+^{\lambda D}}{e_+^{\lambda}} \right)^t =T(D)^t. $\\
2. $ T(QDQ^{-1})=  \frac{e_+^{\lambda Q DQ^{-1}}}{e_+^{\lambda}} = \frac{Qe_+^{\lambda D}Q^{-1}}{e_+^{\lambda}} = T(D). $\\
3. $ T(D_1 \oplus D_2)= \frac{e_+^{\lambda (D_1 \oplus D_2)}}{e_+^{\lambda}} = \frac{e_+^{\lambda D_1} \oplus \ e_+^{\lambda D_2}}{e_+^{\lambda}} =
\frac{e_+^{\lambda D_1 }}{e_+^{\lambda}} \oplus \frac{e_+^{\lambda D_1 }}{e_+^{\lambda}} = T(D_1)\oplus T(D_2).$\\
4. Let $c_j$ be the linear functional on matrices such that $c_j(S)$ is the sum of the entries of the $j$-column of the matrix $S$, i.e.
$c_j(S)=\sum_{i=1}^nS_{ij}.$ A matrix $S$ is column stochastic if its entries are non-negative and $c_j(S)=1$ for all $j$.
It is easy to check that if $S$ is column  stochastic so is $S^k$ for $k\geq 1$. Assume that $D$ is a column stochastic matrix, then
$$c_j(T(D))= c_j\left(  \frac{e_+^{\lambda D}}{e_+^{\lambda}}  \right)=
\frac{\sum_{k=1}^{\infty}\frac{\lambda^k c_j(D^k)}{k!}}{e_+^{\lambda}}=
\frac{\sum_{k=1}^{\infty}\frac{\lambda^k }{k!}}{e_+^{\lambda}}=\frac{e_+^{\lambda}}{e_+^{\lambda}}=1.$$
5. Assume that $[D_1,D_2]=0$, then we have that
$$T(D_1 + D_2)= \frac{e_+^{\lambda (D_1 + D_2)}}{e_+^{\lambda}} = \frac{e^{\lambda (D_1 + D_2)}- I}{e_+^{\lambda}} =
\frac{e^{\lambda D_1} e^{\lambda D_2}- I}{e_+^{\lambda}} = \frac{(e_+^{\lambda D_1}+I) (e_+^{\lambda D_2}+I)- I}{e_+^{\lambda}}$$
$$= \frac{e_+^{\lambda D_1} e_+^{\lambda D_2}}{e_+^{\lambda}} \ + \ \frac{e_+^{\lambda D_1}}{e_+^{\lambda}} \ + \
\frac{e_+^{\lambda D_2}}{e_+^{\lambda}}= e_+^{\lambda}T(D_1)T(D_2) \ + \ T(D_1) \ + \ T(D_2).$$
6. By definition the expected matrix  $E\widehat{D}$ of the random matrix $\widehat{D}$  is given by
$$E\widehat{D}=\sum_{k=1}^{\infty}\widehat{D}(k)p(k)= \sum_{k=1}^{\infty} D^k \frac{\lambda^k}{e_+^{\lambda}k!}=
\frac{\sum_{k=1}^{\infty} \frac{(\lambda D)^k}{k!} }{e_+^{\lambda}} = \frac{e_+^{\lambda D}}{e_+^{\lambda}}=T(D).$$

\end{proof}

Recall that the Poisson probability $P$ on $\mathbb{N}$ is given by $P(k)= e^{-\lambda}\frac{\lambda^k}{k!}.$ Let again $p$ be the probability on
$\mathbb{N}_+$
given by $p(k)=\frac{\lambda^k}{e_+^{\lambda}k!}$.

\begin{prop}{\em With the above notation we have that:
\begin{enumerate}
\item For $k\geq 1$ we have that $p(k)=\frac{P(k)e^{\lambda}}{e^{\lambda}-1}$.

\item Let $X$ be a random variable with distribution $p$, then $$EX=\frac{ \lambda e^{\lambda} }{ e^{\lambda } - 1 }, \ \ \ \ \ \
EX^2=\frac{ (\lambda^2 + \lambda) e^{\lambda} }{ e^{\lambda } - 1
}, \ \ \ \ \ \  \mbox{and} \ \ \ \ \ \ VX=\frac{ \lambda e^{2\lambda} -  (\lambda^2 +
\lambda)e^{\lambda}}{ (e^{\lambda } - 1 )^2}.$$

\item We have that
$$p\left(\left| X - \frac{ \lambda e^{\lambda} }{ e^{\lambda } - 1 }\right| \geq c \right) \leq \frac{ \lambda e^{2\lambda} -  (\lambda^2 + \lambda)e^{\lambda}}{ c^2(e^{\lambda } - 1)^2 }.$$
\end{enumerate}
}
\end{prop}

\begin{proof}
1. is a trivial calculation, and 2. follows from 1. and the well-known fact that if $X$ is a Poisson random variable then
$$EX=\lambda \ \ \ \  \text{and} \ \ \ \ VX=\lambda^2 + \lambda.$$ 3. Direct consequence of the Chebyschev's inequality.
\end{proof}

In order to provide a combinatorial interpretation for the MICMAC, PageRank, Heat Kernel and PWP
methods we need a few definitions, see for example \cite{diaz}. Let $\mathbb{R}$-set
be the category of $R$-weighted finite sets, i.e. the category whose
objects are pairs $(x, \omega)$ where $x$ is a finite set and $\omega: x
\rightarrow \mathbb{R}$ is a map. A morphism in
$\mathbb{R}$-Set from $(x_1, \omega_1)$ to $(x_2, \omega_2)$ is a
map $\alpha: x_1 \rightarrow x_2$ such that $\omega_1 =
\omega_2 \circ \alpha .$  $\mathbb{R}$-set is a distributive category provided with a natural
valuation map $$| \ \ |:\mathbb{R}\text{-}\mathrm{set} \rightarrow \mathbb{R}$$ given by
$$|x, \omega|= \sum_{i \in x}\omega(i) .$$
Note that the definition above may, sometimes, be applied for some infinite sets as well.
If $e$ is an edge of a directed graph we denote by $se$ and
$te$ the starting point and the endpoint of $e$, respectively. A
path $\gamma$ of length $k$ from a vertex $j$ to a vertex $i$  in
a graph  is a sequence of edges $\gamma= (\gamma_1,....,\gamma_k)$
such that $s\gamma_1=j$,  $t\gamma_i=s\gamma_{i+1}$ for $1
\leq i <k$ and $t\gamma_k=i.$ We let $P(i,j)$ be the set of
all paths from $j$ to $i$, and $P_k(i,j)$ be the set of paths
of length $k$ from $j$ to $i$.\\

We assume that the graph of direct influences has associated
matrix $D$. In our applications we use the $\mathbb{R}$-weighted
sets $(P_k(i,j), \omega)$, $(P_k(i,j), \rho)$, $(P(i,j), \sigma)$  and  $(P(i,j),
\tau)$ where the weights $\omega$, $\rho$, $\sigma$ and $\tau$
 are given, respectively, on a path $\gamma=(\gamma_1,
....,
\gamma_k)$ in the graph of direct influences by

\begin{eqnarray*}
\omega (\gamma) &=& \prod_{i=1}^{k}D_{t\gamma_i, s\gamma_i}, \\
\rho (\gamma) &=& \prod_{i=1}^{k}\left( p\overline{D} + (1-p)E_n \right)_{t\gamma_i, s\gamma_i}, \\
\sigma(\gamma) &=& \sigma_k(\gamma)\frac{\lambda^k}{k!} =  \prod_{i=1}^{k}(D-I)_{t\gamma_i, s\gamma_i} \frac{\lambda^k}{k!}.\\
\tau(\gamma) &=& \tau_k(\gamma)\frac{\lambda^k}{e_+^{\lambda}k!} = \left( \prod_{i=1}^{k}D_{t\gamma_i, s\gamma_i}\right) \frac{\lambda^k}{e_+^{\lambda}k!}.
\end{eqnarray*}

The following result provides combinatorial interpretation for the
MICMAC, PageRank, Heat Kernel and the PWP methods.

\begin{thm}\label{ja}{\em
\begin{enumerate}
\item Let $T$ be the MICMAC matrix of indirect influences. Then $T_{ij}=|P_k(i,j), \omega|$.

\item Let $T$ be the PageRank matrix of indirect influences. Then $T_{ij}=\lim_{k \to \infty}|P_k(i,j),\rho|$.

\item Let $T$ is the Heat Kernel matrix of indirect influences.  Then $$T_{ij}=
|P(i,j), \sigma|=\sum_{k=0}^{\infty}|P_k(i,j), \sigma_k| \frac{\lambda^k}{k!}.$$

\item Let $T$ is the PWP matrix of indirect influences. Then $$T_{ij}=
|P(i,j), \tau|=\sum_{k=1}^{\infty}|P_k(i,j), \tau_k| \frac{\lambda^k}{e_+^{\lambda}k!}.$$
\end{enumerate}}
\end{thm}

\begin{proof}1. The MICMAC matrix $T$ is equal to $D^k$, thus we have that
$$T_{ij}=D^k_{ij}= \sum_{s_1,....,s_{k-1}=1}^nD_{is_1}....D_{s_{k-1}j}= $$ $$\sum_{(\gamma_1,
...,\gamma_k)\in P_k(i,j)}\prod_{i=1}^{k}D_{t\gamma_i, s\gamma_i} = \sum_{\gamma \in P_k(i,j)}\omega(\gamma)=|P_k(i,j), \omega|.$$

\noindent 2. The PageRank matrix $T$ is given $\mbox{lim}_{k\to \infty}\left[p\overline{D} + (1-p)E_n\right]^k$, thus we have that

\begin{eqnarray*}
  T_{ij} &=& \lim_{k \to \infty}\left[p\overline{D} + (1-p)E_n\right]_{ij}^k \\
   &=& \lim_{k \to \infty}
 \sum_{s_1,....,s_{k-1}=1}^n\left[p\overline{D} + (1-p)E_n\right]_{is_1}....\left[p\overline{D} + (1-p)E_n\right]_{s_{k-1}j} \\
   &=&  \lim_{k \to \infty} \sum_{(\gamma_1,
...,\gamma_k)\in P_k(i,j)}\prod_{i=1}^{k}\left[p\overline{D} + (1-p)E_n\right]_{t\gamma_i, s\gamma_i} \\
   &=& \lim_{k \to \infty}\sum_{\gamma \in P_k(i,j)}\rho(\gamma)=\lim_{k \to \infty}|P_k(i,j), \rho|.
\end{eqnarray*}

\noindent 3. The Heat Kernel matrix $T$ of indirect influences is given by $T=e^{\lambda(D-I)},$ thus we have that
$$T_{ij}= e^{\lambda(D-I)}_{ij}= \sum_{k=0}^{\infty} (D-I)^k_{ij}\frac{\lambda^k}{k!}=
 \sum_{k=0}^{\infty} |P_k(i,j), \omega_k|\frac{\lambda^k}{k!}=|P(i,j), \omega|.$$

\noindent 4. The PWP matrix $T$ is given by $T=\frac{e_+^{\lambda D}}{e_+^{\lambda}},$ thus we have that
$$T_{ij}= \frac{e_{+ij}^{\lambda D}}{e_+^{\lambda}}= \sum_{k=1}^{\infty} \frac{ D^k_{ij}\lambda^k}{e_+^{\lambda}k!}=
 \sum_{k=1}^{\infty} |P_k(i,j), \tau_k|\frac{\lambda^k}{e_+^{\lambda}k!}=|P(i,j), \tau|.$$

\end{proof}

\section{Examples}

\begin{exmp}{\em
Let $\mathrm{L}_n$ be the linear  graph with $n$ vertices:  $$ 1 \rightarrow 2 \rightarrow 3 \rightarrow ......... \rightarrow n-1 \rightarrow n.$$ Thus $\mathrm{L}_n$
is the graph  with vertex set $[n]$ and an unique edge from
$j$ to $j+1$ for $j<n$. Figure \ref{1} shows the graph $\mathrm{L}_4$. \\

\begin{figure}
\begin{center}
\includegraphics[scale=0.3]{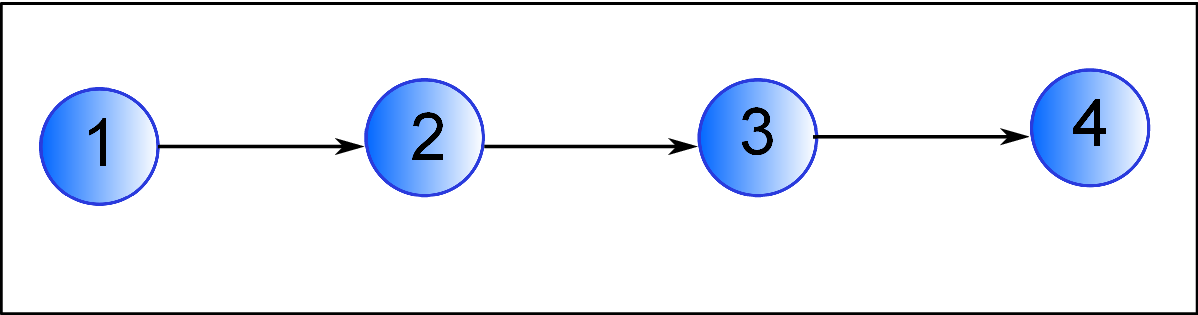}
\caption{Linear graph $\mathrm{L}_4$ with 4 vertices.}
\label{1}
\end{center}
\end{figure}

According to MICMAC, for fixed $k<n$, vertex $j$
will only influence vertex $j+k$. Thus the vertices $\{1,...,n-k
\}$ are the most influential ones each with influence $1$. For $k
\geq n$ the matrix of indirect influences vanishes. For example MICMAC for $\mathrm{L}_3$ and $\mathrm{L}_4$, with $k=4$,
predicts vanishing influences and dependencies. \\

PageRank for $\mathrm{L}_3$
gives the dependency vector $(0.17, 0.34, 0.47)$ making vertex $3$
the most dependent one.  For $n=4$  PageRank
dependencies are $(0.12, 0.21, 0.30, 0.37)$. Again, we see that vertex $4$ is the most dependent one, and that
vertex $1$ has a positive dependency. One can check that according to PageRank  the dependency of vertex
$1$ in the graph $\mathrm{L}_n$ is given by $0.17, \ 0.12, \ 0.08, \ 0.059, \ 0.046 \ $ as the parameter $n$ increases from $3$ to
$7$. For large $n$ we expect PageRank to produce an almost
vanishing dependency for vertex $1$; for small $n$  though this dependency is not quite zero and may be significative. \\

\begin{figure}
\begin{center}
\includegraphics[scale=0.3]{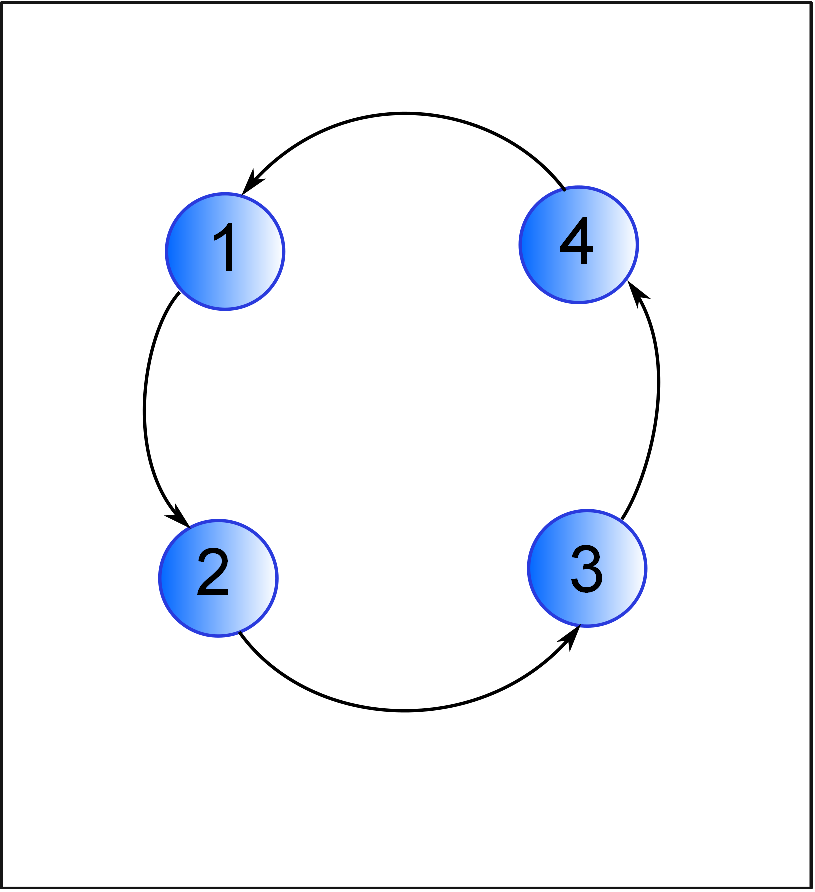}
\caption{Cyclic graph with 4 vertices.}
\label{3}
\end{center}
\end{figure}

For $\mathrm{L}_3$ and $\lambda=0.5$, the Heat Kernet method gives the vectors
$f=(0.98, 0.9, 0.6)$ and $d=(0.6, 0.9, 0.98).$ Notice that vertex $3$ has a non-vanishing influence,
and vertex $1$ a non-vanishing influence. The highest influence and dependency  of a vertex is on itself. Indeed, self-influence
is responsible for more than half of the total influence of a vertex.\\

For $\mathrm{L}_3$ and $\lambda=1$, the PWP method gives the
vectors $f=(1.5, 1, 0)$ and $d=(0,1,1.5).$ Vertex $3$ is the most
dependent and has a vanishing influence; vertex $1$ has a vanishing dependency.
For $\mathrm{L}_4$ and $\lambda=1,$ PWP yields the vectors $f=(\frac{5}{3},
\frac{3}{2}, 1, 0)$ and $d=(0,1,\frac{3}{2},
\frac{5}{3}).$ Thus vertex $4$ is the most dependent one and has vanishing influence; vertex $1$ has a vanishing
dependency and has the highest influence.\\

For arbitrary $n$ and $\lambda$ one can check that PWP  yields the matrix $T$
given by

$$T_{j+s,j}=\left\{
\begin{array}{cc}
\frac{\lambda^s}{e_+^{\lambda}s!} & \mbox{ for \ }1 \leq j <n, \ \ \ 1 \leq s
\leq n-j \\

0 & \mbox{otherwise.}
\end{array}\right.$$
From the expression above we see that for $s\leq n-j-1$ we have
$$\frac{T_{j+s+1,j}}{T_{j+s,j}}=\frac{\lambda}{s+1},$$ and therefore we obtain that
$$T_{j+s+1,j} >T_{j+s,j} \ \ \ \ \ \ \mbox{for} \ \ \ \ s < \lambda -1,$$ $$T_{j+s+1,j}
<T_{j+s,j} \ \ \ \ \ \ \mbox{for} \ \ \ \ s > \lambda -1.$$ Therefore $T_{j+s,j}$  achieves
its maximum, for fixed $j$, at $s=\lfloor \lambda  \rfloor$ if
$\lambda \geq 1$, and at $s=1$ if
$\lambda < 1$.

}
\end{exmp}

\begin{figure}
\begin{center}
\includegraphics[scale=0.3]{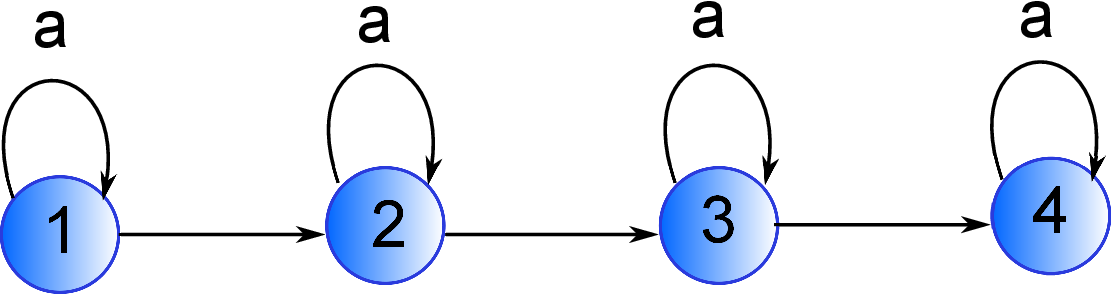}
\caption{Jordan graph with 4 vertices.}
\label{2}
\end{center}
\end{figure}

\begin{exmp}{\em
Let $\mathrm{C}_n$ be the cyclic graph on $[n]$, i.e. the graph with an
edge from $j$ to $j+1$ for $j<n$ and an edge from $n$ to $1$. Figure \ref{3} shows the cyclic graph
$\mathrm{C}_4$. In this case
the four methods yield the same vector of dependencies, but they do so for
quite different reasons. \\

MICMAC makes the influence and dependence of each vertex equal to
$1$. Indeed for fixed $k$,   vertex $j$ will only influence the vertex  $j+k$
mod $n$. \\

The PageRank vector of dependencies is also $(1,....,1)$, indeed
each vertex is equally dependent on every other vertex, i.e. the PageRank matrix $T$ has all its
entries equal to $\frac{1}{n}.$ \\

The Heat Kernel, for $n=4$ and $\lambda= 0.5$, yield the vector of dependencies $(1,1,1,1)$. The influence
of a vertex on the other vertices decreases as the distance in the cyclic ordering increases. The highest influence of a vertex is
on itself.\\

With PWP the matrix $T$ of indirect influences is given, for $j,s \in [n]$ and $j+s$
taken mod $n$,  by
$$T_{j+s,j}=\frac{1}{e_+^{\lambda}}\sum_{k=0}^{\infty}\frac{\lambda^{nk+s}}{(nk+s)!}.$$ Thus for $\lambda=1$ we have that $$T_{j+s,j}=\frac{1}{e_+^{1}}\sum_{k=0}^{\infty}\frac{1}{(nk+s)!},$$
and therefore we get that $$T_{j+1,j} > T_{j+2,j} > .... >
T_{j+n,j}=T_{j,j}.$$ Thus, the influence and the dependence of
a vertex $j$ are both equal to $1$. Vertex $j$ influences all other vertices;
it has a higher influence over the vertices closer to it in the cyclic ordering.
Its lowest influence is on itself. }
\end{exmp}

\begin{exmp}{\em
From Theorem 1, properties 2 and 3, choosing an appropriated basis one can always reduced the computation of
the PWP matrix $T=T(D)$ of indirect influences to the case where $D$ is a Jordan block of a matrix in Jordan canonical
form.  We let $J_n(a)$ be the Jordan graph associated with a Jordan block with value $a$ on the diagonal.
The Jordan graph $\mathrm{J}_4$ is shown in Figure \ref{2}.\\

Thus, we assume that $D$ is a  matrix such that $D_{jj}=a , \  D_{j+1, j}=1, \  \mbox{and} \
D_{ij}=0 \ \mbox{otherwise.}$ Vertex $j$ exerts a non-vanishing indirect
influence over the vertices $j+s$ with $0 \leq s \leq n -j.$
Notice that at vertex $j$ a path can either stay at $j$ or move to
$j+1$, thus  the MICMAC matrix for fixed $k$ is given by
$$|P_k(j+s,j)|={k
\choose s} a^{k-s}.$$

Therefore the PWP matrix $T$ of indirect
influences is given by
$$T_{j+s,s}=\frac{1}{e_+^{\lambda}}\sum_{k=s}^{\infty}{k \choose s}a^{k-s}\frac{\lambda^k}{k!}=
\frac{ e^{a \lambda} \lambda^s}{(e^{\lambda}-1)s!}.$$

}
\end{exmp}

\begin{figure}
\begin{center}
\includegraphics[scale=0.3]{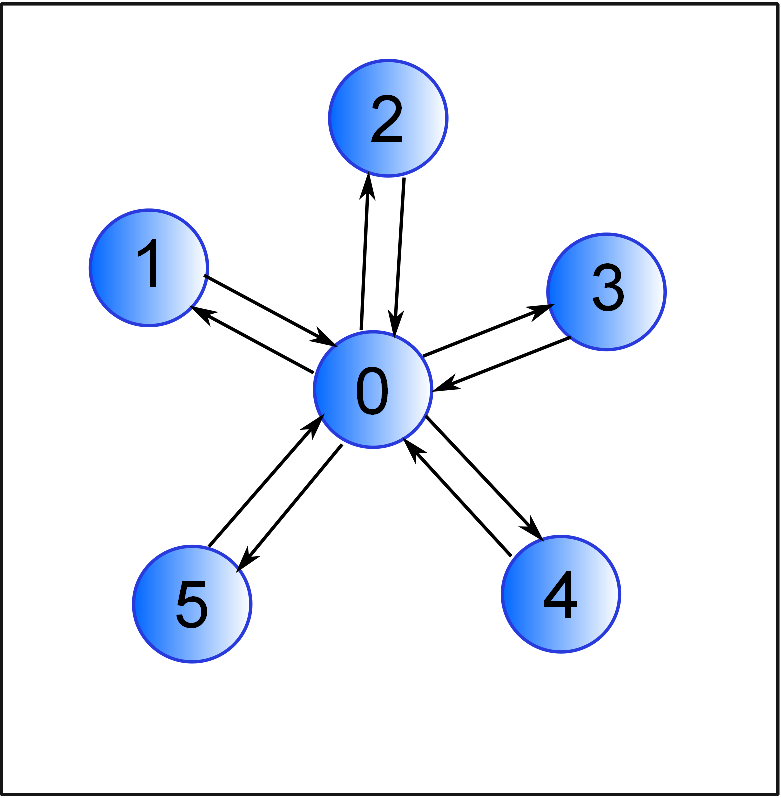}
\caption{Star graph with 6 vertices.}
\label{4}
\end{center}
\end{figure}

\begin{exmp}{\em
Consider the star graph $\mathrm{S}_n$ with vertex set $\{0,1,....,n \}$, see Figure \ref{4}, and directed edges from $0$  to $ i \in [n]$ and viceversa.
It is not hard to see that the PWP matrix of indirect influences $T$ is a symmetric matrix given, for $i,j \in [n],$ by

\begin{eqnarray*}
T_{0,0} &=& \frac{  e^{\lambda \sqrt{n}} + e^{ - \lambda \sqrt{n}   } }{2e_+^{\lambda}}, \\
T_{j,0} &=& \frac{1}{e_+^{\lambda}} \sum_{k=0}^{\infty}n^k\frac{\lambda^{2k+1}}{(2k+1)!}, \\
T_{i,j} &=& \frac{  e^{\lambda \sqrt{n} } + e^{ - \lambda   \sqrt{n} }
}{2ne_+^{\lambda}}.
\end{eqnarray*}

Thus the vertex with greater influence and
dependence is the vertex $0$. PageRank yields a similar result
making $0$ the most dependent vertex.}
\end{exmp}

\section{Final comments}

The main goal of this note is to propose an alternative method for counting indirect influences.
It seems convenient to have a pool of options, as well as a comparative study of the various possibilities.
We introduced the PWP method for counting indirect influences. Applications of PWP to real-world networks is currently underway.
We worked with a discrete scenario where influences are transmitted linearly.
Lifting those restrictions will conduce to continuous non-linear models. This more general setting will be considered elsewhere.

\bigskip

\bigskip

\noindent ragadiaz@gmail.com\\
\noindent Escuela de Matem\'aticas, Universidad Sergio Arboleda, Bogot\'a, Colombia\\

\end{document}